\documentclass[reqno,12pt]{amsart}
\begin{document}
\emergencystretch10pt
\def\ov#1{\overline{#1}}
\def\1#1{{\widetilde{#1}}}
\def\2#1{{\overline{#1}}}
\def\3#1{{\mathcal{#1}}}
\def\d{{\partial}}
\def\reg{{\rm reg}}
\def\id{{\text{\rm id}}}
\def\supp{{\text{\rm supp}}}
\def\eps{\varepsilon}
\def\R{{\mathbb R}}
\def\I{{\Bbb I}}
\def\C{{\mathbb C}}
\def\Z{{\Bbb Z}}
\newtheorem{Th}{Theorem}[section]
\newtheorem{Df}[Th]{Definition}
\newtheorem{Cr}[Th]{Corollary}
\newtheorem{Pp}[Th]{Proposition}
\newtheorem{Lm}[Th]{Lemma}

\title[Domains of holomorphy with edges]{Domains of holomorphy with edges and lower dimensional boundary singularities}

\author[D. Zaitsev and G. Zampieri]{Dmitri Zaitsev and Giuseppe Zampieri}
\address{Mathematisches Institut, Eberhard-Karls-Universit\"at T\"ubingen,
        72076 T\"ubingen, GERMANY.
        E-mail address: dmitri.zaitsev@uni-tuebingen.de}
\address{Dipartimento di Matematica Pura ed Applicata, 
Universit\`a degli Studi di Padova, via G. Belzoni 7, 35131 Padova, ITALY.
        E-mail address: zampieri@math.unipd.it}
\begin{abstract}
Necessary and sufficient geometric conditions are given for domains
with regular boundary points and edges to be domains of holomorphy
provided the remainder boundary subset
is of zero Hausdorff $1$-codimensional measure.
\end{abstract}

\maketitle

\section{Introduction}
The positive solution to the classical Levi problem due 
to {\sc Oka} \cite{O42,O53}, {\sc Bremermann} \cite{B54} and {\sc Norguet} \cite{N54}
asserts that a domain $\Omega\subset\C^N$ whose boundary is of class $C^2$
is a domain of holomorphy provided the Levi form of the boundary 
is everywhere positively semidefinite (see e.g. surveys \cite{S84,P94}).
In contrast to this, for domains with singularities on the boundary
there seems to be a lack of such geometric conditions in the literature.

\subsection{Piecewise smooth domains}
A natural generalization of smooth domains is given by the class of
so-called {\em piecewise smooth} domains
whose boundaries are pieces of hypersurfaces satisfying
suitable transversality conditions
(see e.g. \cite{SH81,P82,N88,K92,F93,MP94}).
However, for the domains of holomorphy,
the class of piecewise smooth domains seems to be very restrictive.
For instance, an envelope of holomorphy of a domain with real-analytic (even algebraic) boundary
does not need to be piecewise smooth as the example 
$$\Omega := \{(z,w)\in\C^2 : |w|^2+(|z|^2-1)^2<2 \}$$
shows. Indeed, the Cauchy formula argument implies that
the envelope of holomorphy (and also convex, polynomially convex and rationally convex hulls)
of $\Omega$ is the union $\Omega\cup\{|z|<1,|w|<\sqrt{2}\}$.

In this paper we consider a larger class of domains $\Omega$,
whose smooth boundary pieces may not be extended
to closed smooth hypersurfaces in a neighborhood of $\d\Omega$.
We also allow singular subsets in the boundary 
that are only controlled to have zero $1$-codimensional 
(with respect to the dimension of the boundary)
Hausdorff measure.

\subsection{The class $L^{2,\infty}$.}
For an open subset $U\subset\R^m$ ($m\ge 1$), 
denote by $L^{2,\infty}(U)$ the space of real continuous functions $h$ on $U$
that are twice continuously differentiable 
with bounded derivatives on an open dense subset of $U$.

We say that a domain $\Omega\subset\R^n$ 
is {\em of class $L^{2,\infty}$} if for every $a\in\d\Omega$
there exists a system of local ($C^2$-smooth) coordinates 
$(x_1,\ldots,x_n)=(x',x_n)\in \R^{n-1}\times\R$ in a
neighborhood $U=U'\times I$ of $a$ and
a function $h\in L^{2,\infty}(U')$ such that
\begin{equation}\label{h}
\Omega\cap U = \{ (x',x_n)\in U : x_n > h(x')\}.
\end{equation}

It is easy to see that every domain with piecewise $C^2$-smooth
boundary is of class $L^{2,\infty}$.

\subsection{Regular and edge points.}
Given a subset $A\subset\R^n$,
we say that a point $a\in A$ is ($C^2$-){\em regular} 
if $A\cap U_a$ is a smooth hypersurface of class $C^2$
for some neighborhood $U_a\subset\R^n$ of $a$. 
If $a\in A$ is not regular,
we call it a ($C^1$-){\em edge point} if
there exists a neighborhood $U_a$ of $a$
and a connected closed $(n-2)$-dimensional submanifold 
$M_a\subset \d\Omega\cap U_a$ of class $C^1$, referred to as an {\em edge at $a$},
that contains all nonregular points of $A\cap U_a$.

It $\Omega\subset\C^N$ is a domain of holomorphy,
it is a standard fact that
the Levi form (see \S\ref{levipar})
at every regular point is positively semidefinite.
The classical example of Hartogs 
$$\Omega := \{|z|<1,|w|<1/2\}\cup \{1/2<|z|< 1,|w|< 1\} \subset\C^2 $$
shows that the converse does not hold in general 
even for piecewise smooth domains.
The edge boundary points of $\Omega$ in this example, 
where all holomorphic functions extend,
are precisely those whose tangent cones are not convex.
Recall that the {\em tangent cone} (in the sense of Whitney)
of $\Omega\subset\C^N$ at a point $a\in\d\Omega$, denoted by $T_a\Omega$,
is defined 
to be the set of all possible limits of $t_k(a_k-a)\in\C^N$,
where $a_k\in \Omega$ and $t_k\in \R_+$
are sequences with $a_k\to a$ as $k\to\infty$.

\subsection{Main results.}
It turns out that, together with the Levi form condition 
for regular points of $\d\Omega$,
the cone convexity for edges points
guarantees that $\Omega$ is a domain of holomorphy.
No condition on the other points is required
provided the set of those points if 
of Hausdorff $(2N-2)$-dimensional measure zero.
More precisely, we have the following result.

\begin{Th}\label{main0}
Let $\Omega\subset\C^N$ be a domain of class $L^{2,\infty}$
and $E\subset\d\Omega$ be a closed 
subset of Hausdorff $(2N-2)$-dimensional measure zero.
Suppose that the following hold:
\begin{enumerate}
\item[(i)] every point of $\d\Omega\setminus E$ 
is either regular or an edge point;
\item[(ii)] at every regular point of $\d\Omega\setminus E$ 
the Levi form is positively semidefinite;
\item[(iii)] at every edge point of $\d\Omega\setminus E$ 
the tangent cone of $\Omega$ is convex.
\end{enumerate}
Then $\Omega$ is a domain of holomorphy.
\end{Th}

By the well-known fact, $\Omega$ is a domain of holomorphy
if and only if the function 
$\psi(z):=-\log d(z,\d\Omega)$
 is plurisubharmonic in $\Omega$,
where $d$ denotes the euclidean distance (see e.g. Theorems 2.6.5 and 4.2.8 in \cite{H}).
In particular, it follows from Theorem~\ref{main0} for $N\ge 2$ that the tangent cone $C_a$ of $\Omega$
at any $a\in \d\Omega$ cannot be strictly concave 
(i.e. the interior of $C_a$ cannot contain a hyperplane).
Indeed, otherwise the function $\psi$ would be equal to $-\log \|z-a\|$ in an open subset of $\Omega$
and hence would not be plurisubharmonic.
This shows, on the other hand, that $\psi$ cannot be directly used to prove Theorem~\ref{main0}
because in Theorem~\ref{main0} there is no convexity condition on the cone at 
the points from the ``exceptional'' subset $E\subset\d\Omega$.
In fact we prove the plurisubharmonicity of the function
\begin{equation}\label{phi}
\phi(z):=-\log(x_{2N}-h(x')) + \lambda\|z\|^2
\end{equation}
near the boundary for some $\lambda>0$
rather than of $\psi$, where $(x_1,\ldots,x_{2N})$ and $h$ satisfy (\ref{h}).

The necessity of the convexity condition (iii) depends on the complex geometry of the edges. We show:

\begin{Pp}\label{real-edge}
Let $\Omega\subset\C^N$ be a domain of holomorphy and 
suppose that for an edge point $a\in\d\Omega$,
there exists an edge $M_a$ which is not a complex hypersurface
in any neighborhood of $a$.
Then the tangent cone of $\Omega$ at $a$ is convex.
\end{Pp}

On the other hand, if an edge can be chosen to be a complex hypersurface, the convexity condition (iii)
does not need to hold as the example of $\Omega:=D\times\C\subset\C^2$ shows
with $D\subset\C$ a nonconvex polygon. 
Therefore we have to distinguish between edge points satisfying the assumptions of 
Proposition~\ref{real-edge} that we call {\em real edge points} and
other edge points $a\in\d\Omega$, where any edge must be locally a complex hypersurface.
In the second case $a$ is said to be a {\em complex edge point}.
Then we impose the convexity condition only at {\em real} edge points.
In this more general situation the above function $\phi$ given by (\ref{phi})
is not always plurisubharmonic. Nevertheless,
we obtain the following necessary and sufficient
geometric conditions for domains to be domains of holomorphy
as a consequence of Theorem~\ref{main0}
and Proposition~\ref{real-edge}.

\begin{Cr}\label{maincor}
Let $\Omega\subset\C^N$ be a domain of class $L^{2,\infty}$
and $E\subset\d\Omega$ be a closed subset of Hausdorff $(2N-2)$-dimensional measure zero.
Suppose that the following hold:
\begin{enumerate}
\item[(i)] every point in $\d\Omega\setminus E$ is either regular or an edge point;
\item[(ii)] for every $a\in \d\Omega$ 
there exist a neighborhood $U_a$ and a complex
hypersurface $N\subset \d\Omega$
that contains all complex edge points in $\d\Omega\cap U_a$.
\end{enumerate}
Then $\Omega$ is a domain of holomorphy if and only if
the Levi form at every regular point $a\in \d\Omega\setminus E$
is positively semidefinite and the tangent cone
at every real edge point $a\in \d\Omega\setminus E$ is convex.
\end{Cr}

Finally we would like to mention that 
the statements of Theorem~\ref{main0} (and of Corollary~\ref{maincor})
also hold for relatively compact domains in Stein manifolds.
Indeed, in this case Theorem~\ref{main0} implies that the domain is 
{\em locally Stein}.
Hence it is a domain of holomorphy by 
a result of {\sc Fornaess} and {\sc Narasimhan} (\cite{FN80}, Theorem 3.1.1).

\section{The Levi form and plurisubharmonicity}\label{levipar}

Recall that the {\em Levi form}
 at a point $a$ of a real function $\rho$ of class $C^2$
in an open subset of $\C^N$
is the Hermitian form defined in local holomorphic coordinates $z=(z_1,\ldots,z_N)$ by
$$L\rho(\xi,\eta)=L\rho(a)(\xi,\eta):=\sum_{k,l}\frac{\d^2\rho}{\d z^k\d\2{z^l}}(a)\,\xi^k \2{\eta^l},
\quad \xi,\eta\in\C^N.$$
We write 
$$\d\rho(a)^\perp:=\{\xi\in \C^N :\d\rho(a)(\xi)=0 \}.$$

The {\em Levi form of a domain} $\Omega$ at a regular point $a\subset\d\Omega$
is the restriction $L\rho|_{\d\rho^\perp}$,
where $d\rho\ne 0$ and $\Omega$ is locally given by $\rho<0$.
The norms are defined in the standard way:
$$\|\d\rho(a)\|:=\sup_{\|\xi\|=1} |\d\rho(a)(\xi)|, \quad 
\|L\rho(a)\|:=\sup_{\|\xi\|=\|\eta\|=1} |L\rho(a)(\xi,\eta)|.$$

\begin{Lm}\label{rho}
Let $\rho<0$ be a negative function of class $C^2$ in an open subset of $\C^N$ 
and $\lambda>0$ be a constant such that the following holds:
\begin{enumerate}
\item[(i)] $L\rho|_{\d\rho^\perp}$ is positive definite;
\item[(ii)] $\|L\rho\|^2\le \lambda(\|\partial\rho\|^2+\rho\|L\rho\|)$.
\end{enumerate}
Then the function $\phi(z):=-\log(-\rho(z)) + \lambda\|z\|^2$ is plurisubharmonic.
\end{Lm}

\begin{proof}
We have
\begin{equation}\label{llog}
L\phi(\xi,\eta)=\rho^{-2}\d\rho(\xi)\2{\d\rho(\eta)} - \rho^{-1}L\rho(\xi,\eta) + 
\lambda \langle \xi,\eta \rangle,
\end{equation}
where $\langle \xi,\eta \rangle :=\xi_1\2\eta_1+\cdots+\xi_N\2\eta_N$.
Every vector $\zeta\in \C^N$ can be written as $\zeta=\zeta_1+\zeta_2$ with
\begin{equation}\label{split}
|\d\rho(a)(\zeta_1)|=\|\d\rho(a)\|\cdot\|\zeta_1\| \text{ and } \d\rho(a)(\zeta_2)=0.
\end{equation}
Applying (\ref{llog}) to $\xi=\eta=\alpha_1\zeta_1+\alpha_2\zeta_2$
with $\zeta_1,\zeta_2$ satisfying (\ref{split}) we obtain
$$L\phi(a)(\alpha_1\zeta_1+\alpha_2\zeta_2,\alpha_1\zeta_1+\alpha_2\zeta_2) = 
(\alpha_1,\alpha_2) (A + B)
\begin{pmatrix}
\2\alpha_1 \\
\2\alpha_2
\end{pmatrix},
$$
where 
$$A=
\begin{pmatrix}
\rho^{-2}\|\d\rho\|^2\cdot\|\zeta_1\|^2 - \rho^{-1} L\rho(\zeta_1,\zeta_1) &  - \rho^{-1} L\rho(\zeta_1,\zeta_2) \\
- \rho^{-1} L\rho(\zeta_2,\zeta_1) & \lambda \|\zeta_2\|^2
\end{pmatrix}
$$
and
$$B=
\begin{pmatrix}
\lambda\|\zeta_1\|^2 & 0 \\
0 & - \rho^{-1} L\rho(\zeta_2,\zeta_2)
\end{pmatrix}.
$$
The matrix $B$ is positively semidefinite by (i). 
It is sufficient to show that $A$ is also positively semidefinite,
i.e. $\det A\ge 0$ by Sylvester's criterion.
But this follows from (ii):
$$\det A \ge \rho^{-2}(\lambda \|\d\rho\|^2 + \lambda \rho\|L\rho\| - \|L\rho\|^2)
\|\zeta_1\|^2\|\zeta_2\|^2\ge 0.$$
\end{proof}

\section{Piecewise plurisubharmonicity}
\subsection{One-dimensional case.}
In the following let $I\subset\R$ denote an open interval and $A\subset I$ a finite subset.
If $f$ is continuously differentiable with bounded derivative on $I\setminus A$,
then for every $a\in S$ there exist one-sided limits
$$f(a)_-:=\lim_{x\to a,\, x<a}f(x) \quad \text{ and } \quad f(a)_+:=\lim_{x\to a,\, x>a}f(x).$$
By elementary calculus we have

\begin{Lm}\label{fsum}
Let $A\subset I$ and $f$ be as before 
and suppose that $f$ has compact support in $I$.
Then 
$$\int_I f' dx + \sum_{a\in A}(f(a)_+-f(a)_-) = 0.$$
\end{Lm}

\begin{Cr}\label{phisum}
Let $A\subset I$ be as before, $\phi\in C^0(I)\cap C^2(I\setminus A)\cap  L^{2,\infty}(I)$
be arbitrary and $\alpha\in C^2(I)$ have compact support in $I$. Then
\begin{equation}\label{phijump}
\int_I (\phi\alpha''-\phi''\alpha)dx = \sum_{a\in A}(\phi'(a)_+-\phi'(a)_-)\alpha(a).
\end{equation}
\end{Cr}

The corollary is obtained by applying Lemma~\ref{fsum} to $f:=\phi\alpha'-\phi'\alpha$.

\subsection{Higher-dimensional case.}\label{higher}
Now consider an open subset $\Omega\subset\R^n$
and let $S\subset\Omega$ be a (locally closed) hypersurface of class $C^1$.
Given a point $a_0\in S$ we fix a neighborhood $U\subset\Omega$ of $a_0$ such that $U\setminus S$ 
has exactly two connected components $U_+$ and $U_-$.

\begin{Lm}\label{one-sided}
Let $\phi\in C^2(U\setminus S)\cap L^{2,\infty}(U)$ be a function in $U$.
Then for every $v\in \R^n$ the directional derivatives 
$D_v\phi|_{U_-}$ and  $D_v\phi|_{U_+}$ extend Lipschitz-continuously
to $U_-\cup S$ and $U_+\cup S$ respectively with the one-sided limits
\begin{equation}\label{limits}
D_v\phi(a)_-:=\lim_{x\to a,\, x\in U_-} D_v\phi(x) \quad \text{ and } \quad 
D_v\phi(a)_+:=\lim_{x\to a,\, x\in U_+} D_v\phi(x)
\end{equation}
for $a\in S$.
Moreover, if $\phi$ is in addition continuous on $U$, one has
$D_v\phi(a)_-=D_v\phi(a)_+$ whenever $v$ is tangent to $S$ at $a$.
In particular, the sign of the expression
\begin{equation}\label{jump}
D_v\phi(a)_+-D_v\phi(a)_-
\end{equation}
is independent of the choice of a transversal vector $v$ pointing into $U_+$.
\end{Lm}

\begin{proof}
Boundedness of the first and second derivatives of $\phi$ on $U_-$ and $U_+$
implies the one-sided Lipschitz extendibility of $\phi$ and its first derivatives.
In particular, if $\phi$ is continuous on $U$,
then the restriction $\phi|_S$ coincides with both one-sided limits.
Hence for $v$ tangent to $S$, one has $D_v\phi_+=D_v(\phi|_{S})= D_v\phi_-$ 
as required.
\end{proof}

We observe that, if we interchange $U_-$ with $U_+$, the sign of (\ref{jump}) remains the same,
because also $v$ (pointing into $U_+$) changes the sign. 
This consideration motivates the following definition.

\begin{Df}\label{jdef}
A function $\phi\in C^2(U\setminus S)\cap L^{2,\infty}(U)$, continuous in $U$,
is said to be {\em transversally convex} at $S$
if the expression {\rm (\ref{jump})} is non-negative for any $a\in S$ and any $v$ pointing into $U_+$.
\end{Df}

In the following we write $\3H^m$ for the Hausdorff $m$-dimensional measure.

\begin{Pp}\label{derivatives}
Let $\Omega\subset\R^n$ be open, $G$ be closed in $\Omega$ with $\3H^{n-1}(G)=0$,
$S\subset\Omega$ be a hypersurface 
of class $C^1$ with $\2S\setminus S\subset G$ and  
$$\phi\in C^0(\Omega)\cap C^2(\Omega\setminus(S\cup G))\cap L^{2,\infty}(\Omega)$$
be a function which is transversally convex at $S$.
Then for every non-negative function $\alpha\in C^2(\Omega)$
with compact support in $\Omega$ the inequality
\begin{equation}\label{basicin}
\int_\Omega \phi\, \frac{\d^2\alpha}{\d x_j^2}\, dx \ge 
\int_\Omega \frac{\d^2\phi}{\d x_j^2}\, \alpha\, dx
\end{equation}
holds for every $j=1,\ldots,n$.
\end{Pp}

{\bf Remark.} Since all functions in (\ref{basicin}) are measurable and bounded,
both integrals exist with respect to the Lebesgue measure on $\R^n$.

\begin{proof} Without loss of generality we may assume $j=n$ and, by 
taking a suitable partition of unity,
$${\rm supp}\, (\alpha) \subset U'\times I \subset \Omega$$
for an open subset $U'\subset\R^{n-1}$ and an interval $I\subset\R$.
We write $x=(x',x_n)\in \R^{n-1}\times\R$.
Since $\3H^{n-1}(G)=0$, we see from Fubini's theorem that 
$G\cap (\{x'\}\times I) = \emptyset$
for all $x'$ outside a zero measure subset $G'\subset U'$.
Furthermore, by Sard's theorem applied to the projection $S\mapsto \R^{n-1}$, the vector $v:=(0,1)$
is not tangent to $S$ at the points of $S\cap (\{x'\}\times I)$ for all $x'$
outside a zero measure subset $G''\subset U'$.
Then for $x'\notin G'\cup G''$ the set 
$$A_{x'}:=\{x_n\in I : (x',x_n)\in (S\cap{\rm\supp}(\alpha)) \}$$ 
is finite.
Thus we can apply Corollary~\ref{phisum} to the restriction of $\phi$ and $\alpha$ to 
$\{x'\}\times I$. According to Definition~\ref{jdef}, 
the right-hand side in (\ref{phijump}) is non-negative.
Then the required inequality is obtained by integrating over $U'\setminus (G'\cup G'')$
the non-negative left-hand side of (\ref{phijump}).
\end{proof}

\begin{Cr}\label{cor}
Under the assumptions of Proposition~{\rm \ref{derivatives}} suppose that $\Omega\subset\C^N$
and $\phi$ is plurisubharmonic in $\Omega\setminus(S\cup G)$. 
Then $\phi$ is plurisubharmonic in the whole $\Omega$.
\end{Cr}

\begin{proof}
Given a vector $\xi\in \C^N$ we can find linear complex coordinates $z_k=x_k+iy_k$ ($k=1,\ldots,N$)
such that $\xi=(0,\ldots,0,1)$. Then
\begin{equation}\label{leviin}
L\phi(a)(\xi,\xi)=\frac{\d^2\phi}{\d z_N\d\2{z_N}}(a) \xi^N\2{\xi^N} = 
\frac{1}{2} \left(\frac{\d^2\phi}{\d x_N^2}(a) + \frac{\d^2\phi}{\d y_N^2}(a) \right) \ge 0
\end{equation}
for all $a\in \Omega\setminus(S\cup G)$ by the plurisubharmonicity of $\phi$ there.
Since $S\cup G\subset\Omega$ is of zero measure, (\ref{leviin}) implies
\begin{equation}\label{basicsum0}
\int_\Omega \frac{\d^2\phi}{\d x_N^2}\, \alpha\, dz \2{dz}
+ \int_\Omega \frac{\d^2\phi}{\d y_N^2}\, \alpha\, dz \2{dz} \ge 0
\end{equation}
for every non-negative function $\alpha\in C^2(\Omega)$ with compact support in $\Omega$.
By Proposition~{\rm \ref{derivatives}},
\begin{equation}\label{basicsum}
2\int_\Omega \phi\: L\alpha(a)(\xi,\xi) dz \2{dz} =
\int_\Omega \phi\, \frac{\d^2\alpha}{\d x_N^2}\, d\Omega 
+ \int_\Omega \phi\, \frac{\d^2\alpha}{\d y_N^2}\, d\Omega  \ge 0.
\end{equation}
Since $\xi$ is arbitrary, we conclude that $\phi$ has a non-negative Levi form in distributional sense. 
By the continuity of $\phi$,
the last fact is equivalent to the plurisubharmonicity of $\phi$ on the whole $\Omega$
(see \cite{H}, Theorem 1.6.11).
\end{proof}

\section{Proof of Theorem~\ref{main0}}\label{mainproof}
Suppose that conditions (i) and (ii) hold.
Since a domain in $\C^N$ is a domains of holomorphy if it is pseudoconvex
(see \cite{H}, Theorem 4.2.8)
and due to the local characterization of pseudoconvexity (\cite{H}, Theorem 2.6.10)
it is sufficient to show that every point $x_0\in\d\Omega$ has a neighborhood $U$
such that $\Omega\cap U$ is pseudoconvex.
We use the standard identification $\C^N\cong \R^n$ with $n:=2N$.
Since $\Omega$ is of class $L^{2,\infty}$, we can choose a neighborhood
$U=U'\times I\subset\R^{n-1}\times\R$ of $x_0$ and a continuous function
$h\colon U'\to I$ of class $L^{2,\infty}$ such that (\ref{h}) holds.
Define a continuous function $\rho\colon U\to \R$ 
with $\Omega\cap U=\{\rho<0\}$ by
$$\rho(x',x_n):=h(x')-x_n.$$

Set $\1\Omega:=\Omega\cap U$ and denote by 
$\1\Omega_\reg$ the subset of all $x=(x',x_n)\in \1\Omega$
such that $(x',h(x'))\in \d\Omega$ is ($C^2$-)regular.
We wish to apply Lemma~\ref{rho} to the restriction of $\rho$ to $\1\Omega_\reg$.
To check the assumption (i) in Lemma~\ref{rho} we observe that
\begin{equation}
L\rho(x',x_n)=L\rho(x',h(x')), \quad \d\rho(x',x_n)=\d\rho(x',h(x'))
\end{equation}
for all $x\in \1\Omega_\reg$. 
Hence it follows by condition (ii) in Theorem~\ref{main0}
that $L\rho(x)|_{\d\rho(x)^\perp}$ is positively semidefinite
whenever $x\in \1\Omega_\reg$. In order to satisfy the assumption (ii) in Lemma~\ref{rho}
we shrink the neighborhood $U$ of $x_0$
such that 
\begin{equation}\label{rhoeq}
|\rho(x)| \ \|L\rho(x)\|\le 1/8
\end{equation}
holds for all $x\in \1\Omega_\reg$.
This is possible because the second derivatives of $\rho$ are bounded on $\1\Omega_\reg$
 and $\rho$ is continuous with $\rho(x_0)=0$.
Since $\|\d\rho(x)\|\ge 1/2$ for $x\in \1\Omega_\reg$,
(\ref{rhoeq}) implies 
$$\|\d\rho\|^2+\rho\|L\rho\|\ge 1/4 - 1/8 = 1/8$$
on $\1\Omega_\reg$ and hence the existence of $\lambda>0$
such that the assumption (ii) in Lemma~\ref{rho} in satisfied.
By Lemma~\ref{rho}, the function 
$$\phi\colon \1\Omega\to \R, \quad \phi(z):=-\log(-\rho(z))+\lambda \|z\|$$
is plurisubharmonic in $\1\Omega_\reg$.

Next we wish to apply Corollary~\ref{cor} to $\phi$ in $\1\Omega$.
For this we first construct a $C^1$ hypersurface $S\subset\1\Omega$
satisfying the assumptions. By using the definition, 
for every edge point $a\in\d\Omega\cap U$,
we can choose a neighborhood $U_a=U'_a\times I\subset\R^{n-1}\times\R$ 
with $U'_a\subset \R^{n-1}$ a euclidean ball
and a closed $(n-2)$-dimensional real submanifold $M_a\subset \d\Omega\cap U_a$ of class $C^1$ 
such that all nonregular points of 
$\d\Omega\cap U_a$ are contained in $M_a$.
Furthermore we can shrink $U_a$ and $M_a$ such that all nonregular points of 
$\d\Omega\cap \2{U_a}$ are contained in $\2{M_a}$ and
$\3H^{n-2}(\2{M_a}\setminus M_a)=0$.
By the choice of the coordinates $(x',x_n)$, the projection
$$M'_a:=\{x'\in\R^{n-1} : \exists x_n, (x',x_n)\in M_a\}$$
is a closed submanifold of $U'_a$.
In view of condition (i) any nonregular point of $\d\Omega\setminus E$
is an edge point.
By considering a sequence of all rational points in $\C^N$ 
we may choose a sequence $a_m$ ($1\le m <\infty$) of edge points 
such that the union $\cup_m (V'_{a_m}\times I)$ covers all nonregular points in 
$(\d\Omega\setminus E)\cap U$,
where $V'_{a_m}\times I$ is an open neighborhood of $a_m$ and $\2{V'_{a_m}}\subset U'_{a_m}$.
Then by applying Fubini's theorem, we may shrink the balls $U'_{a_m}$
and their submanifolds $M'_{a_m}$
to obtain the additional property 
\begin{equation}\label{addprop}
\3H^{n-2}(\d U'_{a_m}\cap M'_{a_k})=0 \text{ for all } m,k\ge 1.
\end{equation}
Define $U'_m:=\cup_{j\le m} U'_{a_j}$.

We claim that there exist increasing sequences of $C^1$ submanifolds $S'_m\subset U'_m$
and closed subsets $E'_m\subset U'_m$ with $\3H^{n-2}(E'_m)=0$ such that 
\begin{enumerate}
\item[(i)] $S'_k\cap U'_m=S'_m$ and $E'_k\cap U'_m=E'_m$ for all $k\ge m$,
\item[(ii)] all nonregular points in $(\d\Omega\setminus E)\cap (\2{U'_m}\times I)$ outside $E'_m\times I$
 are contained in $\2{S'_m}\times I$,
\item[(iii)] $\3H^{n-2}(\2{S'_m}\setminus S'_m)=0$.
\end{enumerate}
We construct the sequences $S'_m$ and $E'_m$ by induction on $m$.
For this we set $S'_1:=M'_{a_1}$, $E'_1:=\emptyset$. 
If $S'_{m-1}$ and $E'_{m-1}$ are already constructed, define
$$S'_m:=S'_{m-1}\cup (M'_{a_m}\cap (U'_{a_m}\setminus \2{U'_{m-1}}))$$
and $E'_m:=E_{m-1}\cup (\2{M'_{a_m}}\cap\d U'_{m-1})$.
It is easy to see from our construction that (i)-(iii) are satisfied.

Define 
$$S:=((\cup_m S'_m)\times I)\cap \1\Omega  \text{ \ and \ }
F:=((\cup_m E'_m)\times I)\cap\1\Omega.$$
Then $\3H^{n-1}(F)=0$ and $F$ is closed in $\1\Omega$ by (i).
It follows from (i) that $S$ is a $C^1$ hypersurface in $\1\Omega$.
Furthermore, $\3H^{n-1}(\2{S}\setminus S)=0$ by (iii).

Denote by $E'\subset U'$ the projection of the subset $E\subset\d\Omega$.
Since $\d\Omega\cap U$ is the graph of $h$, $E'$ is closed in $U'$.
Next we subtract $E'\times I$ from $S$ and denote the remainder again by $S$.
Then for every $x=(x',x_n)\in S$, $(x',h(x'))\in \d\Omega$
is either $C^2$-regular or an edge point. Finally we define 
$$G:= F\cup (\2{S}\setminus S)\cup ((E'\times I)\cap\1\Omega),$$
where the closure of $S$ is taken in $\1\Omega$.
By (ii), $\1\Omega\setminus (S\cup G)\subset \1\Omega_\reg$
and therefore $\phi$ is plurisubharmonic in $\1\Omega\setminus (S\cup G)$.

In order to apply Corollary~\ref{cor} it remains to show that $\phi$
is transversally convex at $S$ (see Definition~\ref{jdef}).
But this is a direct consequence of condition (iii) in Theorem~\ref{main0}
on the convexity of tangent cones.
We conclude that $\phi$ is plurisubharmonic on $\1\Omega$.
Then, for a sufficiently small euclidean ball $B(x_0,\eps)$ centered at $x_0$,
$$\max \big(\phi(z),-\log(\eps-||z-x_0||)\big)$$
is a plurisubharmonic
exhaustion function for $\Omega\cap B(x_0,\eps)$.
This shows pseudoconvexity of $\Omega\cap B(x_0,\eps)$
and hence completes the proof.

\section{Proof of Proposition~\ref{real-edge}}
Let $a\subset\d\Omega$ be a real edge point satisfying the assumptions
of Proposition~\ref{real-edge}.
As before we choose a neighborhood
$U=U'\times I\subset\R^{n-1}\times\R$ of $a$ and a continuous function
$h\colon U'\to I$ of class $L^{2,\infty}$ such that (\ref{h}) holds.
Let $M_a\subset\d\Omega$ be an edge and denote by $M'_a$ its projection on $U'$.
As in \S\ref{higher} choose $U'$ sufficiently small such that 
$M'_a$ divides $U'$ into two parts $U'_-$ and $U'_+$.
Then the tangent cone is given by 
\begin{equation}\label{cone}
T_a\Omega = \{(v',v_n)\in\R^n : v_n\ge D_{v'}h(a)_\pm \text{ for } 
v'\in T_aU'_\pm\},
\end{equation}
where the one-sided limits exist by Lemma~\ref{one-sided}.

We first suppose that $M_a$ is generic at $a$, i.e. $T_aM_a+iT_aM_a=\C^N$.
We prove the statement by contradiction 
assuming that the tangent cone of $\Omega$ at $a$ is not convex.
Then there exists a linear disc 
$A\colon \Delta\to T_a\Omega$, $t\mapsto tv$
with $v\in T_aM_a$ and $A(\d\Delta\setminus\{-1,1\})$ in the interior of $T_a\Omega$.
Here $\Delta:=\{|t|<1\}\subset\C$.
Furthermore, for $\zeta:=(0,1)\in\R^{n-1}\times\R$ and $\eps>0$ sufficiently small,
the ``deformed disc'' $A'(t):=tv + t^2 \eps\zeta$ sends the whole boundary $\d\Delta$
to the interior of $T_a\Omega$. Finally for $\mu>0$ sufficiently small
the ``rescaled disc'' $A''(t):=a+\mu A'(t)$ sends the boundary $\d\Delta$ into $\Omega$.
Then the Cauchy formula argument for $A''$ and its translations shows that all holomorphic functions
in $\Omega$ extend holomorphically to a neighborhood of $a=A''(0)$ which contradicts
the assumption that $\Omega$ is a domain of holomorphy.

Now consider the general case.
Let $M_a$ be an edge at $a$ satisfying the assumptions.
Then $M_a$ contains generic points arbitrarily close to $a$.
Hence the tangent cones at those points are convex
by the above argument.
But then the explicit formula (\ref{cone}) shows
that the cone $T_a\Omega$ is also convex completing the proof.

\section{Proof of Corollary~\ref{maincor}}
The necessity of the Levi form condition is well-known.
The necessity of the convexity follows from Proposition~\ref{real-edge}.
For the converse it is sufficient to show the local pseudoconvexity
at every point $a\in\d\Omega$ as in \S\ref{mainproof}.
Let $N\subset U_a$ be given by condition (ii).
By a coordinate change we may assume that $N$ is locally given by $z_N=0$.
Then it follows from (\ref{cone}) that,
if $U_a$ is sufficiently small polydisc centered at $a$, 
the intersection $\Omega\cap U_a$ can be mapped via
$$(z_1,\ldots,z_{N-1},z_N)\mapsto (z_1,\ldots,z_{N-1},z_N^\alpha)$$
biholomorphically onto
a domain satisfying assumptions of Theorem~\ref{main0}.
Here $\alpha$ is a sufficiently small positive number.
The required conclusion follows now from Theorem~\ref{main0}.

\end{document}